\documentclass{article}

\usepackage{amssymb, amsmath}
\usepackage[matrix, arrow, curve]{xy}
\newcommand \lap {\lambda^{\prime}}
\newcommand \pip {\pi^{\prime}}

\newcommand \la {\lambda}
\newcommand \gz {{\cal G}}

\newcommand \R {{\cal R}}
\newcommand \A {{\cal A}}

\newcommand \q {{\bf q}}
\newcommand \Prob {\mathbb P}
 \newcommand \U {{\cal U}}

\newcommand \Y {{\cal Y}}

\newcommand \F {{\cal F}}

\newcommand \wab {{\cal W}_{{\cal A},B}}

\newcommand  {\modk} {{\cal M}_{\kappa}}

\newcommand \wper   {{\cal W}_{per}}

\newcommand \pper   {{\cal P}_{per}}

\newcommand \wpi {{{\cal W}}_{per}^-(\pi)}

\newcommand \wpip {{{\cal W}}_{per}^+(\pi)}

\newcommand \omabz {{\Omega}_{{\cal A},B}^{\mathbb Z}}

\newcommand \omab {{\Omega}_{{\cal A},B}}

\newcommand \vone {{\cal V}^{(1)}}

\newtheorem {lemma} {Lemma}

\newtheorem {theorem} {Theorem}

\newtheorem {proposition} {Proposition}
\newtheorem{corollary}{Corollary}
\begin{document}
\title{Logarithmic asymptotics for the number of 
periodic orbits of the Teichm{\"u}ller flow 
on Veech's space of zippered rectangles.}
\author{Alexander I. Bufetov.}
\date{}
\maketitle

\section{Introduction}

The aim of this note is to obtain a logarithmic asymptotics for 
number of periodic orbits on Veech's space of zippered rectangles, such that the norm of the 
corresponding renormalization matrix does not exceed a given value.
The space of zippered rectangles is a finite branched cover 
of a set of full measure in the moduli space of abelian differentials with prescribed 
singularities. In particular, all periodic orbit in moduli space can be lifted to the space of 
zippered rectangles. The renormalization matrix corresponds to the action of 
the corresponding pseudo-Anosov automorphism in relative homology of the underlying surface with respect 
to the singularities of the abelian differential.

Let $\R$ be a Rauzy class of permutations on $m$ symbols, 
let $\Omega_0(\R)$ be Veech's moduli space of zippered rectangles,
and let $P^t$ be the Teichm{\"u}ller flow on $\Omega_0(\R)$. 
As shown by Veech, to each periodic orbit $\gamma\subset \Omega_0(\R)$ 
there corresponds a renormalization matrix $A(\gamma)\in SL(m, {\mathbb Z})$
(see the next section for precise definitions). 
The period of the orbit $\gamma$
is the logarithm of the spectral radius of $A(\gamma)$. 
The aim of this note is to give an asymptotics of the number of orbits 
such that {\it the norm} of the renormalization matrix does not exceed $\exp(T)$.
More precisely, for a matrix $A$, write 
$$
||A||=\max_{j}\sum_i |A_{ij}|,
$$
and denote by $Per(\R,T)$ the set of  periodic orbits $\gamma$ for 
$P^t$ such that  $||A(\gamma)||\leq\exp(T)$. 

\begin{theorem}
\label{main}
Let $\R$ be a Rauzy class of permutations on $m$ symbols. Then
$$
 \lim_{T\to\infty}\frac {\log \# Per(\R, T)}{T} = m.
$$
\end{theorem}

{\bf Remark.} By a Theorem of Veech \cite{veechteich}, 
the entropy of the flow $P^t$ is equal to $m$.

The proof proceeds as follows. Veech's moduli space of zippered 
rectangles $\Omega_0(\R)$ admits a natural Lebesgue measure class, and the 
Teichm{\"u}ller flow $P^t$ preserves an absolutely continuous
finite ergodic invariant measure $\mu_{\R}$ on $\Omega_0(\R)$.

The flow $P^t$  has a pair of 
complementary ``stable'' and ``unstable'' foliations. As Maryam Mirzakhani has pointed out to me, 
the measure $\mu_{\R}$ has the Margulis property \cite{margulis} 
of uniform expansion on unstable leaves. This observation allows to obtain the 
logarithmic asymptotics of the number of periodic orbits whose period does 
not exceed $T$ ``in the compact part'' of the 
space $\Omega_0(\R)$, using mixing of the flow $P^t$, due to Veech \cite{veechteich},
and following the classical argument of Margulis \cite{margulis}. Note that ``in the compact part'' 
of the moduli space the period of the orbit is comparable (up to an additive constant) to the logarithm of the norm 
of the renormalization matrix.

To estimate the number of periodic orbits ``at infinity'', it is convenient to represent 
the Teichm{\"u}ller flow as a suspension flow over the natural extension 
of the Rauzy-Veech-Zorich induction map. Each periodic 
orbit is then coded by a finite word over a countable alphabet.
A compact set in the moduli space can be chosen in such a way that periodic 
orbits passing through it correspond to words that contain a given subword.
But then, to a word coding a periodic orbit, one can assign its concatenation with 
the word, corresponding to the compact set. The norm of the corresponding matrix 
grows at most by a mutltiplicative constant under this procedure. 
The asymptotics of periodic orbits in the whole space is thus reduced to 
the asymptotics of orbits passing ``through the compact part''.

By a Theorem of Veech \cite{veech}, the moduli space $\Omega_0(\R)$ of zippered rectangles
admits an almost everywhere defined, surjective up to a set of measure zero, 
finite-to-one projection map on the moduli space of abelian 
differentials with prescribed singularities 
(the genus and the orders of the singularities are uniquely determined by the Rauzy class 
$\R$ \cite{kontsevichzorich}). 
The Teichm{\"u}ller flow $P^t$ projects to the Teichm{\"u}ller flow 
on moduli space of abelian differentials, and periodic orbits are taken to 
periodic Teichm{\"u}ller geodesics (perhaps of smaller period). 
It would be interesting to see if a similar asymptotics could be obtained 
in moduli space of abelian differentials as well.

\section{Rauzy-Veech-Zorich induction.}

The Teichm{\"u}ller flow on Veech's space of zippered rectangles
can be represented as a suspension flow over the natural extension of the
Rauzy-Veech-Zorich induction map on the space of interval
exchange transformations.
Since the Rauzy-Veech-Zorich induction map has a
natural countable Markov partition,
we obtain a symbolic coding for the Teichm{\"u}ller
flow \cite{veech, zorich, bufetov}.
In particular, periodic orbits are represented by periodic
symbolic sequences.
This representation will be essential for our argument.

In this section, we briefly recall the definitions of the Rauzy-Veech-Zorich 
induction, of Veech's space of zippered rectangles, and of the symbolic representation
for the Teichm{\"u}ller flow. For a more detailed presentation, see \cite{veech, zorich, bufetov}.

\subsection{Rauzy operations $a$ and $b$.}

Let $\pi$ be a permutation on $m$ symbols. The permutation $\pi$ will always be assumed
irreducible in that $\pi\{1,\dots ,k\}=\{1,\dots,k\}$ iff $k=m$.

Rauzy operations $a$ and $b$ are defined by the formulas:

$$
a\pi(j)=\begin{cases}
\pi j,&\text{if $j\leq \pi^{-1}m$;}\\
                 \pi m,&\text{if $j=\pi^{-1}m+1$;}\\
                  \pi(j-1),&\text{ other $j$.}
\end{cases}
$$

$$
b\pi(j)=\begin{cases}
\pi j,&\text{if $\pi j\leq \pi m$;}\\
                 \pi j+1,&\text{if $\pi m<\pi j<m$;}\\
                  \pi m+1,&\text{ if $\pi j=m$.}
\end{cases}
$$

These operations preserve irreducibility.
The {\it Rauzy class} of a permutation
$\pi$ is defined as the set of all permutations that can be obtained
from $\pi$ by repeated application of the operations $a$ and $b$.

For $i,j=1,\dots,m$,
denote by $E_{ij}$ an $m\times m$ matrix of which
the $i,j$-th element is equal to $1$, all others to $0$.
Let $E$ be the $m\times m$-identity matrix.

Introduce the matrices
$$
A(a, \pi)=\sum_{i=1}^{\pi^{-1}(m)}E_{ii}+E_{m, \pi^{-1}m+1}+
\sum_{i=\pi^{-1}m+1}^m E_{i,i+1},
$$
$$
A(b, \pi)=E+E_{m,\pi^{-1}m}
$$

For a vector $\la\in{\mathbb R}^m$, $\la=(\la_1, \dots, \la_m)$,
we write
$$
|\la|=\sum_i \la_i.
$$
Let $\Delta_{m-1}$ be the unit simplex in ${\mathbb R}^m$:
$$
\Delta_{m-1}=\{\la\in {\mathbb R}_+^m: |\la|=1\}.
$$
{\it The space} $\Delta(\R)$ {\it
of interval exchange transformations}, corresponding
to a Rauzy class $\R$, is defined by the formula
$$
\Delta(\R)=\Delta_{m-1}\times \R.
$$
Denote
$$
\Delta_{\pi}^+=\{\la\in\Delta_{m-1}| \ \la_{\pi^{-1}m}>\la_m\},\ \ 
\Delta_{\pi}^-=\{\la\in\Delta_{m-1}| \ \la_m>\la_{\pi^{-1}m}\},
$$
and
$$
\Delta^+ = \cup_{\pip\in{\cal R}(\pi)}\Delta_{\pip}^+,\ \ 
\Delta^- = \cup_{\pip\in{\cal R}(\pi)}\Delta_{\pip}^-.
$$

{\it The Rauzy-Veech induction} is a map
$$
{\cal T}: \Delta(\R)\to \Delta(\R),
$$
defined by the formula
$$
{\cal T}(\la,\pi)=\begin{cases}
(\frac{A(\pi, a)^{-1}\la}{|A(\pi, a)^{-1}\la|}, a\pi), &\text{if $\la\in\Delta^-$;}\\
(\frac{A(\pi, b)^{-1}\la}{|A(\pi, b)^{-1}\la|}, b\pi), &\text{if $\la\in\Delta^+$.}
\end{cases}
$$
Veech \cite{veech} showed that the Rauzy-Veech induction
has an absolutely continuous ergodic invariant measure on $\Delta(\R)$; that measure is,
however, infinite.

Following Zorich \cite{zorich}, define the function $n(\la,\pi)$ in the following way.
$$
n(\la, \pi)=\begin{cases}
\inf \{k>0:{\cal T}^k(\la,\pi)\in\Delta^-\},&\text{if $\la\in\Delta_{\pi}^+$;}\\
  \inf \{k>0: {\cal T}^k(\la,\pi)\in\Delta^+\},&\text{if $\la\in\Delta_{\pi}^-$.}
\end{cases}
$$

{The Rauzy-Veech-Zorich induction} is defined by the formula
$$
{\cal G}(\la,\pi)={\cal T}^{n(\la,\pi)}(\la,\pi).
$$

\begin{theorem}[Zorich\cite{zorich}]
The map ${\cal G}$ has an ergodic invariant probability measure, absolutely
continuous with respect to the Lebesgue measure class on $\Delta(\R)$.
\end{theorem}

This invariant measure will be denoted by $\nu$.

\subsection{Symbolic dynamics for the induction map.}

This subsection describes, following \cite{bufetov}, 
the symbolic dynamics for the map $\gz$.

Consider the alphabet
$$
{\cal A}=\{ (c,n,\pi) | \  c=a \ {\rm or}\  b, n\in {\mathbb N}\}.
$$

For $w_1\in\A$, $w_1=(c_1, n_1, \pi_1)$, we write $c_1=c(w_1), \pi_1=\pi(w_1), n_1=n(w_1)$.
For $w_1,w_2\in\A$, $w_1=(c_1, n_1, \pi_1)$, $w_2=(c_2, n_2, \pi_2)$, define
the function $B(w_1, w_2)$ in the following way:
$B(w_1, w_2)=1$ if $c_1^{n_1}\pi_1=\pi_2$ and $c_1\neq c_2$ and $B(w_1, w_2)=0$ otherwise.

Introduce the space of words
$$
{\cal W}_{{\cal A}, B}=\{w=w_1\dots w_n | \
w_i\in{\cal A}, B(w_i, w_{i+1})=1 \  {\rm for \ all}\  i=1, \dots,n\}.
$$

For a word $w\in\wab$, we denote by $|w|$ its length, i.e., the number
of symbols in it;
given two words $w(1), w(2)\in\wab$, we denote by $w(1)w(2)$ their
concatenation.
Note that the word $w(1)w(2)$  need not
belong to $\wab$, unless a compatibility condition is satisfied by
the last symbol of $w(1)$ and the first symbol of $w(2)$.

To each word assign the corresponding renormalization matrix as follows. 
For $w_1\in\A$, $w_1=(c_1, n_1, \pi_1)$, set 
$$
A(w_1)=A(c_1, \pi_1)A(c_1, c_1\pi_1)\dots A(c_1, c_1^{n_1-1}\pi_1),
$$
and for $w\in W_{\A,B}$, $w=w_1\dots w_n$, set 
$$
A(w)=A(w_1)\dots A(w_n).
$$

Words from $\wab$ act on permutations from $\R$: namely, 
if $w_1\in \A$, $w_1=(c_1, n_1, \pi_1)$, then we set $w_1\pi_1=c_1^{n_1}\pi_1$. 
For permutations $\pi\neq \pi_1$, the symbol $w_1\pi$ is not defined. 
Furthermore, for 
$w\in W_{\A,B}$, $w=w_1\dots w_n$, we write
$$
w\pi=w_n(w_{n-1}(\dots w_1\pi)\dots ),
$$
assuming the left-hand side of the expression is defined. 
Finally, if $\pip=w\pi$, then we also write $\pi=w^{-1}\pip$.

Say that $w_1\in\A$ is compatible with $(\la,\pi)\in\Delta(\R)$ if 

\begin{enumerate}
\item either $\la\in\Delta_{\pi}^{+}$, $c_1=a$, and $a^{n_1}\pi_1=\pi$
\item or $\la\in\Delta_{\pi}^{-}$, $c_1=b$, and $b^{n_1}\pi_1=\pi$.
\end{enumerate}

Say that a word $w\in W_{\A,B}$, $w=w_1\dots w_n$ is compatible with $(\la,\pi)$ if 
$w_n$ is compatible with $(\la,\pi)$.  We shall also sometimes say that $(\la,\pi)$ is 
compatible with $w$ instead of saying that $w$ is compatible with $(\la,\pi)$.
We can write 
$$
\gz^{-n}(\la,\pi)=
\{t_w(\la,\pi): |w|=n \ {\rm and}\  w\  {\rm is\  compatible \ with\ } \ (\la,\pi)\}.
$$
Recall that the set $\gz^{-n}(\la,\pi)$ is infinite.

Now consider the sequence spaces 
$$
\Omega_{\A, B}=\{\omega=\omega_1\dots \omega_n\dots |\  \omega_n\in\A, \ B(\omega_n, \omega_{n+1})=1  \ {\rm for \ all}\  n\in {\mathbb N}\},
$$
and
$$
\Omega_{\A, B}^{\mathbb Z}=\{\omega=\dots \omega_{-n}\dots \omega_1\dots \omega_n\dots |\  \omega_n\in\A, \ 
B(\omega_n, \omega_{n+1})=1
\ {\rm for \ all}\  n\in{\mathbb Z}\}.
$$

Denote by $\sigma$ the right shift on both these spaces.
By a Theorem of Veech \cite{veech}, 
the dynamical systems $(\omab, \sigma, \Prob)$ and $\Delta(\R), \gz, \nu)$ are isomorphic. 
Indeed, let $w_1\in\A$ and define the set $\Delta(w_1)$ in the following way.
If $w_1=(a, n_1, \pi_1)$, then
$$
\Delta(w_1)=\{(\la,\pi)\in \Delta^- | \
\exists \lap\in\Delta_{a^{n_1}\pi}^+ {\rm \ such \ that} \ \la=\frac{A(w_1)\lap}{|A(w_1)\lap|}\}.
$$
If $w_1=(b, n_1, \pi_1)$, then
$$
\Delta(w_1)=\{(\la,\pi)\in \Delta^+ | \
\exists \lap\in\Delta_{b^{n_1}\pi}^- {\rm \ such \ that} \ \la=\frac{A(w_1)\lap}{|A(w_1)\lap|}\}.
$$

In other words, for a letter $w_1=(c_1, n_1, \pi_1)$,
the set $\Delta(w_1)$ is the set of all interval exchanges $(\la,\pi)$
such that $\pi=\pi_1$
such that the application of the map $\gz$ to $(\la,\pi)$
results in $n_1$ applications of the Rauzy operation $c_1$.

The coding map $\Phi:\Delta(\R)\to\Omega_{\A, B}$ is given by the formula
\begin{equation}
\Phi(\la,\pi)=\omega_1 \dots \omega_n \dots \ {\rm if} \ 
\gz^n(\la,\pi)\in\Delta(\omega_n).
\end{equation}

The $\gz$-invariant smooth probability measure $\nu$ projects under
$\Phi$ to a $\sigma$-invariant measure
on $\Omega_{\A, B}$; probability with respect to
that measure will be denoted by $\Prob$.

For $w\in W_{\A,B}$, $w=w_1\dots w_n$, let
$$
C(w)=\{\omega\in\Omega_{\A,B} | \ \omega_1=w_1, \dots, \omega_n=w_n\}.
$$
and
$$
\Delta(w)=\Phi^{-1}(C(w)).
$$

W. Veech \cite{veech} has proved the following
\begin{proposition}
Let $\q\in\wab$ be a word such that all entries of the
matrix $A(\q)$ are positive.
Let $\omega\in\omab$ be a sequence having infinitely many occurrences of the word $\q$.
Then the set $\Phi^{-1}(\omega)$ consists of one point.
\end{proposition}

We thus obtain an almost surely bijective symbolic dynamics for the map $\gz$.

\begin{lemma}
\label{gmargulis}
Let $\q$ be a word such that all entries of the matrix $A(\q)$ are positive.
Then there exists a constant $\alpha(\q)$, depending only on $\q$, such that
for any word $w\in\wab$ of the form $w=\q{\tilde w}\q$, ${\tilde w}\in\wab$,
we have
$$
\alpha(\q)^{-1}\leq \Prob(C(w))||A(w)||^m\leq \alpha(\q).
$$
\end{lemma}

Proof. Recall \cite{veech} that the Lebesgue measure of the set $\Delta(w)=\Phi^{-1}(C(w))$
is given by the expression
$$
\frac 1{(\Pi_j\sum_i A(w)_{ij})}.
$$
Since $A(w)=A(\q)A({\tilde w})A(\q)$, there exists
a constant $\beta(\q)$, depending only on $\q$, such that
$$
\beta(\q)^{-1}\leq \frac {||A(w)||^m}{(\Pi_j\sum_i A(w)_{ij})}\leq \beta(\q).
$$

Finally, since the density of the invariant measure for the
Rauzy-Veech-Zorich map is bounded
both from above and from below on $\Delta(\q)=
\Phi^{-1}(C(\q))$ (see \cite{zorich}, \cite{bufetov}),
there exists a constant $\alpha(\q)$, depending only on $\q$, such that
$$
\alpha(\q)^{-1}\leq \Prob(C(w))||A(w)||^m\leq \alpha(\q),
$$
and the Lemma is proved.

\subsection{Veech's space of zippered rectangles.}

A {\it zippered rectangle} associated to the Rauzy class $\R$ is a quadruple $(\la, h, a, \pi)$,  where $\la \in 
{\mathbb R}_+^m, 
h\in{\mathbb R}^m_+, a\in{\mathbb R}^m, \pi\in{\cal R}$, and the vectors 
$h$ and $a$ satisfy the following equations and inequalities (one introduces 
auxiliary components $a_0=h_0=a_{m+1}=h_{m+1}=0$, and sets $\pi(0)=0$, 
$\pi^{-1}(m+1)=m+1$): 

\begin{equation}
\label{zipone}
h_i-a_i=h_{\pi^{-1}(\pi(i)+1)}-a_{\pi^{-1}(\pi(i)+1)-1}, i=0, \dots, m
\end{equation}
\begin{equation}
\label{ziptwo}
h_i\geq 0,   i=1, \dots,m, \  
a_i\geq 0,  i=1, \dots, m-1,  
\end{equation}
\begin{equation}
\label{zipthree}
a_i\leq \min (h_i, h_{i+1}) {\rm \ for\ } i\neq m,i\neq \pi^{-1}m,  
\end{equation}
\begin{equation}
\label{zipfour}
a_m\leq h_m, \ a_m\geq -h_{\pi^{-1}m},\ 
a_{\pi^{-1}m}\leq h_{\pi^{-1}m+1} 
\end{equation}

The area of a zippered rectangle is given by the expression 
$$
\la_1h_1+\dots+\la_mh_m.
$$
Following Veech, we denote by $\Omega(\R)$ the space of all zippered 
rectangles, corresponding to a given Rauzy class $\R$ and satisfying the condition 
$$
\la_1h_1+\dots+\la_mh_m=1.
$$
We shall denote 
by $x$ an individual zippered rectangle.

Veech \cite{veech} defines a flow $P^t$ and a map $\U$ on the space of zippered
rectangles by the formulas:
$$
P^t(\la,h,a,\pi)=(e^t\la, e^{-t}h, e^{-t}a, \pi).
$$
$$
\U(\la,h,a,\pi)=\begin{cases} (A^{-1}(a,\pi)\la, A^t(a,\pi)h, a^{\prime}, a\pi), &\text{if $(\la,\pi)\in\Delta^-$}\\
 (A^{-1}(b,\pi)\la, A^t(b,\pi)h, a^{\prime\prime}, b\pi),&\text{if $(\la,\pi)\in\Delta^+$},
\end{cases}
$$
where
$$
a^{\prime}_i=\begin{cases} a_i, &\text{if $j<\pi^{-1}m$,}\\
h_{\pi^{-1}m}+a_{m-1}, &\text {if $i=\pi^{-1}m$,}\\
a_{i-1}, &\text{other $i$}.
\end{cases}
$$
$$
a^{\prime\prime}_i=
\begin{cases} a_i, &\text{if $j<m$,}\\
-h_{\pi^{-1}m}+a_{\pi^{-1}m-1}, &\text {if $i=m$.}
\end{cases}
$$

The map $\U$ is invertible; $\U$ and $P^t$ commute (\cite{veech}).
Denote
$$
\tau(\lambda, \pi)=\log(|\la|-\min(\la_m,\la_{\pi^{-1}m})),
$$
and for $x\in \Omega(\R)$, $x=(\la, h, a ,\pi)$, write
$$
\tau(x)=\tau(\la,\pi).
$$
Now, following Veech \cite{veech}, define
$$
\Y(\R)=\{x\in\Omega(\R)\; |\; |\la|=1\}.
$$
and
$$
\Omega_0(\R)=\bigcup_{x\in\Y(\R), 0\leq t< \tau(x)}P^tx.
$$
$\Omega_0(\R)$ is a fundamental domain for $\U$ and,
identifying the points $x$ and $\U x$ in $\Omega_0(\R)$, we obtain
a natural flow, also denoted by $P^t$, on $\Omega_0(\R)$.
The space $\Omega_0(\R)$ will be referred to as {\it Veech's moduli space of zippered rectangles}, 
and the flow $P^t$ as the Teichm{\"u}ller flow on the space of zippered rectangles.

The space $\Omega(\R)$ has a natural Lebesgue measure class and so does
the transversal $\Y(\R)$.
Veech \cite{veech} has proved the following Theorem.

\begin{theorem}[Veech \cite{veech}]
\label{zipmeas}
There exists a measure $\mu_{\R}$ on $\Omega(\R)$, absolutely continuous with respect to Lebesgue, preserved by both the map
$\U$  and the flow $P^t$ and such that $\mu_{\R}(\Omega_0(\R))<\infty$.
\end{theorem}

The construction of this measure is recalled in the Appendix.

\subsection{Symbolic representation for the flow $P^t$}

Following Zorich \cite{zorich},  define
$$
\Omega^+(\R)=\{ x=(\la, h,a,\pi)| \ (\la,\pi)\in\Delta^+, a_m\geq 0 \}.
$$
$$
\Omega^-(\R)=\{ x=(\la, h,a,\pi)| \ (\la,\pi)\in\Delta^-, a_m\leq 0 \},
$$
$$
\Y^+(\R)=\Y(\R)\cap \Omega^+(\R),\  \Y^-(\R)=\Y(\R)\cap \Omega^-(\R),\  \Y^{\pm}(\R)=\Y^+(\R)\cup \Y^-(\R).
$$
Take $x\in\Y^{\pm}(\R)$,  $x=(\la,h,a,\pi)$, and let
$\F(x)$ be the first return map of the flow $P^t$ on the transversal $\Y^{\pm}(\R)$.
The map $\F$ is a lift of the map $\gz$ to the space of zippered rectangles: 
\begin{equation}
{\rm if }\ \F(\la,h,a,\pi)=(\lap, h^{\prime}, a^{\prime}, \pip), \ {\rm then} \  (\lap, \pip)=\gz(\lap, \pip).
\end{equation}

Note that if $x\in\Y^+$, 
then $\F(x)\in \Y^-$, and if $x\in\Y^-$, then $\F(x)\in
\Y^+$).

The map $\F$ preserves a natural absolutely continuous invariant measure on
$\Y^{\pm}(\R)$: indeed,
since $\Y^{\pm}(\R)$ is a transversal to the flow $P^t$,
the measure $\mu_{\R}$ induces an absolutely
continuous measure ${\overline \nu}$ on $\Y^{\pm}(\R)$;
since $\mu_{\R}$ is both $\U$
and $P^t$-invariant, the measure ${\overline \nu}$ is $\F$-invariant.
Zorich \cite{zorich} proved that the measure ${\overline \nu}$ is finite
and ergodic for $\F$.

The dynamical system $(\Y^{\pm}, {\overline \nu}, \F)$ is
measurably isomorphic to the system $(\omabz, \Prob, \sigma)$
\cite{veech, bufetov}.
The Teichm{\"u}ller flow $P^t$ on the space
$\Omega_0(\R)$ is a suspension flow over the map $\F$ on $\Y^{\pm}$.
Identifying $\Y^{\pm}$ and $\omabz$, we obtain a symbolic dynamics
for the Teichm{\"u}ller flow in the space of zippered rectangles.

\subsection{Periodic orbits in Veech's space
of zippered rectangles.}

A word $w\in \wab$ will be called {\it admissible} if $ww\in\wab$ and
if for some $r\in{\mathbb N}$ all entries of the matrix
$A(w)^r$ are positive. Note that under these conditions the sequence
$$
{\overline \omega}(w)=\dots w\dots w\dots w\dots
$$
belongs to $\omabz$ and there is a unique interval exchange
transformation, corresponding to it.
We denote by $\wper$ the set of all admissible words.

From the definition of the Teichm{\"u}ller flow on the space of zippered rectanles, we 
immediately have the following 
\begin{proposition}
\begin{enumerate}
\item Let $\gamma$ be a periodic orbit for the flow $P^t$.
Then, the intersection $\gamma\cap\Y^{\pm}$  is not empty and every
$x\in\gamma\cap\Y^{\pm}$ is a periodic point for the map $\F$.
\item
Let $x\in\Y^{\pm}$ be a periodic point for the map $\F$.
Then, there exists an admissible word $w\in\wab$,
such that the symbolic sequence, corresponding to $x$, has the form
$$
\dots w\dots w\dots w \dots
$$
The symbolic sequence
$\dots \dots w\dots w\dots w \dots$ defines the point $x$ uniquely.
\item
If $w\in\wab$ is admissible, then there exists a
unique periodic orbit of the flow $P^t$ on $\Omega_0(\R)$,
corresponding to $w$.
\end{enumerate}
\end{proposition}

The periodic orbit, corresponding to an admissible word $w$
will be denoted by $\gamma(w)$; the length of the orbit $\gamma(w)$
will be denoted by $l(w)$.
As the Proposition above shows, every periodic orbit
$\gamma$ of the flow $P^t$ has the form
$\gamma(w)$ for some $w\in\wper$.
If $\gamma=\gamma(w)$, then we say that $w$ is
a word {\it coding} the periodic orbit $\gamma$.

\section{ Counting admissible words.}

Take $w\in\wper$. Consider the one-sided infinite sequence
$$
\omega(w)=w\dots w\dots\in \omab,
$$
and set
$$
(\la(w),\pi(w))=\Phi^{-1}(\omega(w)).
$$
Periodicity of the sequence $\omega(w)$ implies that $\la(w)$ is an
eigenvector of $A(w)$. Since
all entries of a power of the matrix $A(w)$ are positive,
$\la(w)$ is  an eigenvector
with the maximal eigenvalue of the matrix $A(w)$. 

We denote by $l(w)$ the logarithm of the maximal eigenvalue of the matrix $A$, and we have then
$$
A(w)\la(w)=\exp(l(w))\la(w).
$$
Of course, we have $l(w)\leq \log ||A(w)||$.

For $T>0$, set
$$
\wper(T)=\{w\in\wper| \  \log ||A(w)||\leq T\}.
$$

\begin{lemma}
\label{mainlemma}
The number of admissible words with norm not exceeding $\exp(T)$  
satisfies the following logarithmic asymptotics:
$$
\lim_{T\to\infty}\frac{\log \# \wper(T)}{T}=m.
$$
\end{lemma}

{\bf Remark.} Note that the results of the Lemma 
and of Theorem\ref{main} do not depend on the specific matrix norm used.

Lemma \ref{mainlemma} will be proven in the next section.
Now we derive Theorem \ref{main} from
Lemma \ref{mainlemma}.

This derivation is not automatic, because,
for a given periodic orbit, the coding word $w$ is not unique:
if $w=w_1\dots w_n$ codes a periodic orbit $\gamma$, then
the words $w^{(2)}=w_2\dots w_nw_1$, $w^{(2)}=w_3\dots w_nw_1w_2$,
$\dots$, $w^{(n)}=w_nw_1\dots w_{n-1}$ are precisely all words
coding $\gamma$.

In \cite{bufetov} it is proven that 
the norm of the matrix $A(w)$ grows
exponentially as a function of the number of symbols of the word $w$.
More precisely, Lemma 14 and Corollary 9 in  \cite{bufetov} 
yield the following.
\begin{lemma}
There exists a constant $\alpha_{11}$, depending on the Rauzy class $\R$ only,
such that the following is true. Let $w\in\wab$, $w=w_1\dots w_n$.
Then
$$
||A(w)||\geq \exp(\alpha_{11}n).
$$
\end{lemma}

One can therefore estimate from above 
the number of words coding the same
periodic orbit as a power of the length of the orbit.

\begin{corollary}
There exist  constants $C_{21}, \alpha_{21}$, 
depending on the Rauzy class $\R$ only,
such that the following is true.
Let $\gamma$ be a periodic orbit of length $T$ of the Teichm{\"u}ller
flow on the space of zippered rectangles. Let $w$ be a coding word
for $\gamma$.
Then the number of symbols in $w$ does not exceed $C_{21}T^{\alpha_{21}}$.
\end{corollary}
\begin{corollary}
There exists constants $C_{31}, \alpha_{31}$, depending on the Rauzy class $\R$ only,
such that the following is true.
Let $\gamma$ be a periodic orbit of length $T$ of the Teichm{\"u}ller
flow on the space of zippered rectangles.
Then the number of words, coding $\gamma$, does not 
exceed $C_{31}T^{\alpha_{31}}$.
\end{corollary}

In view of these Corollaries, 
Lemma \ref{mainlemma} implies Theorem \ref{main}.
We now proceed to the proof of Lemma \ref{mainlemma}.

\section{Proof of Lemma \ref{mainlemma}.}

\subsection{`` The compact part'': Margulis's argument.}

Take a word $\q\in\wper$, $\q=q_1\dots q_{2l+1}$,
such that all entries of the matrix
$A(q_1\dots q_l)$ are positive and all the entries of the matrix
$A(q_{l+1}\dots q_{2l+1})$ are positive.
We shall first count the asymptotics of the number of
periodic orbits whose coding words contain the subword $\q$.

We begin with two following simple observations.

\begin{proposition}
Let $w\in\wper$  have the form 
$w={\bf q}{\tilde w}$ for some ${\tilde w}\in\wab$.
Then $w\in\wper$ if and only if the concatenation 
$w{\bf q}$ also lies in $\wab$.
\end{proposition}

\begin{proposition}
\label{qclosing}
There exists a constant $C(\q)$, depending only on $\q$ and such that 
the following is true.
Let $w\in\wper$  have the form 
$w={\bf q}{\tilde w}$ for some ${\tilde w}\in\wab$. 
Then, for any $(\la,\pi)\in\Delta(\q)$,  
we have 
$$
\frac1 {C(\q)}\leq \frac{|A(w)\la|}{\exp(l(w))}\leq C(\q).
$$
\end{proposition}

Now we note that for such words the period of the corresponding orbit 
is comparable (up to an additive constant) to the norm of the renormalization matrix. 

\begin{corollary}
\label{qpernorm}
Let $w\in \wper$  have the form 
$w={\bf q}{\tilde w}$ for some ${\tilde w}\in\wab$. 
Then there exists a positive constant $c(\q)$, depending on $\q$ only and such 
that  
$$
\frac 1{c(\q)}\leq \frac{\exp(l(w))}{||A(w)||}\leq c(\q).
$$
\end{corollary}

The proof immediately follows from the fact that 
$(\la(w), \pi(w))\in\Delta(\q)$ and therefore a constant 
$C(\q)$, depending only on $\q$, may be found in such a way that for all 
$i,j=1,\dots,m$ and for any $w\in\wab$ of the form 
$w={\bf q}{\tilde w}$ we have
$$
\frac 1{ C(\q)}\leq \frac{\la_i(w)}{\la_j(w)}\leq C(\q).
$$

Set
$$
\wper(\q)=\{w\in\wab| \  \exists w(1)\in\wab \  {\rm such \ that \ }
w=\q w(1)\}.
$$
For $T>0$, set
$$
\wper(\q,T)=\{w\in\wper(\q)\;|\; ||A(w)||\leq T\},
$$
and, for an interval $[r,s]\subset {\mathbb R}$, set
$$
\wper(\q,[r,s])=\{w\in\wper(\q)\;|\; r\leq ||A(w)||\leq s\}.
$$

We will also need the quantities

$$
\pper(\q,T)=\{w\in\wper(\q)\;|\; l(w)\leq T\},
$$
and, for an interval $[r,s]\subset {\mathbb R}$, set
$$
\pper(\q,[r,s])=\{w\in\wper(\q)\;|\; r\leq l(w)\leq s\}.
$$
We have 
$$
\#\wper(\q,T)\leq \#\pper(\q,T),
$$
and, by Corollary \ref{qpernorm}, there exists a constant $\alpha_{35}(\q)$ 
depending on $\q$ and such that 
$$
\#\pper(\q,T)\leq \#\wper(\q, T+\alpha_{35}).
$$

\begin{lemma}
\label{margulis}
There exist positive constants
$ \alpha_{41}, \alpha_{42}, T_0$, depending on
$\q$ only, such that for all $T>T_0$ we have
$$
\alpha_{41}^{-1}\exp(mT)\leq \#\pper(\q,[T-\alpha_{42},T+\alpha_{42}])
\leq \alpha_{41}\exp(mT).
$$
\end{lemma}

\begin{corollary}
\label{marguliscol}
$$
\lim_{T\to\infty}\frac{\log \#\wper(\q,T)}{T}=m.
$$
\end{corollary}

Proof of lemma \ref{margulis}.
As in Margulis's classical argument \cite{margulis}, 
the statement of the Lemma follows
from the mixing for the Teichm{\"u}ller flow and the uniform expansion
property for conditional measures on stable leaves.
The proof below is a symbolic version of that argument, and  the almost sure identification of the
transversal $\Y^{\pm}(\R)$ in $\Omega_0(\R)$ and the 
space $\omabz$ is used. By a slight abuse of notation, 
it is convenient to consider the space $\omabz$ itself embedded
into $\Omega_0(\R)$ and to speak of
cylinders, etc., in $\omabz$,
meaning corresponding subsets in $\Y^{\pm}$.

For a word $w\in\wab$, $|w|=n$, $w=w_1, \dots, w_n$,
and an integer $a\in{\mathbb Z}$, denote
$$
C([a, a+n-1], w)=\{\omega\in \omabz| \
\omega_a=w_1, \dots, \omega_{a+n-1}=w_n\}.
$$
For an interval $(r,s)\subset {\mathbb R}$, set
$$
C([a, a+n-1], w, (r,s))=\bigcup_{t\in(r,s)}P^tC([a,a+n-1], w).
$$

For brevity,  set
$$
C^{\prime}(\q)=C([-l,l], q),
$$
and, for an interval $(r,s)\subset {\mathbb R}$, write 
$$
C^{\prime}(\q, (r,s))=\bigcup_{t\in(r,s)}C^{\prime}(\q).
$$

There exists $\epsilon>0$, depending only on $\q$ and such that
$$
P^t(C^{\prime}(\q))\cap C^{\prime}(\q)=\emptyset.
$$
for all $t\in [-10\epsilon, 10\epsilon]$.
Consider the set
$$
C^{\prime}(\q, [-\epsilon, \epsilon])\subset\Omega_0(\R).
$$

The measure
$\mu_{\R}(C^{\prime}(\q, [-\epsilon, \epsilon]))$ is a constant,
depending only on $\q$.

By a theorem of Veech \cite{veechteich}, the Teichm{\"u}ller flow
is strongly mixing, and the measure
$$
\mu_{\R}((C^{\prime}(\q, [-\epsilon, \epsilon])\cap
P^TC^{\prime}(\q, [-\epsilon, \epsilon]))
$$
tends to a  constant, depending only on $\q$, as $T\to\infty$.
Therefore, the measure 
$$\mu_{\R}((C^{\prime}(\q, [-\epsilon, \epsilon])\cap
\bigcup_{t\in [T-\epsilon, T+\epsilon]}
P^t C^{\prime}(\q, [-\epsilon, \epsilon]))
$$
also tends to a  constant, depending only on $\q$.

Since $\epsilon$ only depends on $\q$, using Proposition \ref{qclosing} and Corollary \ref{qpernorm}, we conclude that there exist constants $\alpha_{51}(\q), \alpha_{52}(\q)$, depending on
$\q$ only and such that
\begin{equation}
\label{compar}
\alpha_{51}(\q)^{-1}\leq \bigcup_ {w\in\wper
(\q, [T-\alpha_{52}(\q),T+\alpha_{52}(\q)]))}\Prob(w\q)\leq \alpha_{51}(\q).
\end{equation}

Now recall that, by Lemma \ref{gmargulis}, there exists
a constant $\alpha_{61}(\q)$, depending only on $\q$ and such that for any
$w\in \wper(\q)$, we have
$$
\alpha_{61}(\q)^{-1}\leq \Prob(w\q)||A(w\q)||^m\leq \alpha_{61}(\q).
$$

{\bf Remark.} The fact that the measure of a cylinder, corresponding
to a periodic orbit, is comparable to a power of the length of the
orbit, is the manifestation, in  symbolic language, of the uniform
expansion property.

Since
$$
||A(w)||\leq ||A(w\q)||\leq ||A(w)|| ||A(\q)||,
$$
if $w\in \wper(\q)$, then a constant
$\alpha_{71}(\q)$, depending only on $\q$, can be chosen
in such a way that we have
$$
\alpha_{71}(\q)^{-1}\leq \Prob(w\q)||A(w)||^m\leq \alpha_{71}(\q),
$$
whence, by Lemma \ref{qpernorm}, there exists a constant $\alpha_{81}(\q)$, depending
only on $\q$, such that for any $w\in\wper(\q,[T-\epsilon,T+\epsilon]) $,
we have
\begin{equation}
\label{compar2}
\alpha_{81}(\q)^{-1}\leq \Prob(w\q)\exp(mT)\leq \alpha_{81}(\q).
\end{equation}

From (\ref{compar}), (\ref{compar2}), it follows
that there exist constants $\alpha_{41}(\q), \alpha_{42}(\q)$,
depending only on $\q$ and such that
$$
\alpha_{41}(\q)^{-1}\leq
\frac{\#\{w\in\pper(\q, [T-\alpha_{42}(\q),T+\alpha_{42}(\q)]))\}} {\exp(mT)}
\leq \alpha_{41}(\q).
$$
The proof of Lemma \ref{margulis} is complete.

\subsection{Periodic orbits ``at infinity''.}

To count periodic orbits ``at infinity'', 
the following trick is used. Take a fixed admissible word $\q$ such that all entries of the matrix $A(\q)$ 
are positive.
Assume that an admissible word $w$ is such that the word $\q w$ is also admissible. 
To the word $w$ assign the word $\q w$. 
Since
$$
||A(\q w)||\leq ||A(\q)||||A(w)||,
$$
the logarithmic asymptotics for $\# \wper(T)$ is the same as 
that for $\#\wper(\q,T)$, already known by Lemma \ref{margulis}.

Note, however, that concatenation of two admisible words may not belong to $\wab$.
One therefore considers words beginning and ending at the same permutation.

For $\pi\in\R$,
let $\wpi$ be the set of all words $w\in\wab$,
$w=w_1\dots w_k$, $w_i=(c_i, n_i, \pi_i)$ such that
$$
c_1=a, c_k=b, \pi_1=\pi, b^{n_k}\pi_k=\pi,
$$
and let $\wpip$ be the set of all words $w\in\wab$,
$w=w_1\dots w_k$, $w_i=(c_i, n_i, \pi_i)$ such that
$$
c_1=b, c_k=a, \pi_1=\pi, a^{n_k}\pi_k=\pi.
$$

By definition,
\begin{equation}
\label{unionpi}
\wper=\bigcup_{\pi\in\R} \wpi \cup \wpip.
\end{equation}

Note, however, that, for any $\pi\in\R$,
$$
{\rm if \ }w(1), w(2)\in\wpi \ {\rm then }\ w(1)w(2)\in \wpi.
$$
and, similarly,
$$
{\rm if \ }w(1), w(2)\in\wpip \ {\rm then }\ w(1)w(2)\in \wpip.
$$

For $\pi\in\R$, choose an arbitrary  $\q_{\pi}\in \wpi$,
$\q_{\pi}=q_1\dots q_{2l+1}$,
such that all entries of the matrix $A(q_1\dots q_l)$, as well as
all entries of the matrix $A(q_{l+1}\dots q_{2l+1})$, are positive.

Denote
$$
\wper(\q_{\pi})=\{ w\in \wpi| \ \exists {\tilde w}
\in \wab: w=\q_{\pi}{\tilde w}  \}.
$$

Note that in this case the word ${\tilde w}$
must also belong to $\wpi$.
For $T>0$, set
$$                                                                         
\wper^-(\pi, T)=\{w\in \wpi| \   ||A(w)||\leq T\},
$$
$$
\wper(\q_{\pi},T)=\{w\in \wper(\q_{\pi})| \   ||A(w)||\leq T\}.
$$

Note that if $w\in \wper(\pi,T)$, then
$\q_{\pi}w\in \wper(\q_{\pi}, T+l(\q_{\pi}))$.

For any $T>0$, we have therefore
\begin{equation}
\label{paragon}
\# \wper(\q_{\pi}, T+l(q_{\pi}))
\geq \# \wper^-(\pi, T)
\geq \# \wper(\q_{\pi}, T).
\end{equation}

By Corollary \ref{marguliscol}, for any $\pi\in \R$ we have
\begin{equation}
\label{limword}
\lim_{T\to\infty}\frac {\log\# \wper(\q_{\pi}, T)}{T}=m.
\end{equation}
Equations (\ref{paragon}) and (\ref{limword}) together imply
$$
\lim_{T\to\infty} \frac {\log\# \wper^-(\pi,T)}{T}=m.
$$

An identical argument can be conducted for words in 
$\wpip$ as well, and for the set 
$\wper^+(\pi, T)=\{w\in \wpip| \   l(w)\leq T\}$
we also have the asymptotics 
$$
\lim_{T\to\infty} \frac {\log\# \wper^+(\pi,T)}{T}=m.
$$
The identity(\ref{unionpi}) finally yields
$$
\lim_{T\to\infty} \frac {\log\# \wper(T)}{T}=m.
$$

Lemma \ref{mainlemma} is proved, and 
Theorem \ref{main} is proved.

\section{Appendix A: Zippered rectangles and abelian differentials.}

Veech \cite{veech} established the following connection
between zippered rectangles and moduli of abelian differentials.
A detailed description of this connection is given in \cite{kontsevichzorich}.

A zippered rectangle naturally defines a Riemann surface 
endowed with a holomorphic differential. This correspondence preserves area.
The orders of the singularities of $\omega$ are uniquely 
defined  by the Rauzy class of the permutation $\pi$ (\cite{veech}).
For any $\R$ we thus have a map 
$$
\pi_{\R}: \Omega(\R)\rightarrow\modk,
$$

where ${\kappa}$ is uniquely defined by $\R$.

Veech \cite{veech} proved

\begin{theorem}[Veech]
\label{zipmodule}
\begin{enumerate}
\item Upto a set of measure zero, 
the set $\pi_{\R}(\Omega_0(\R))$ is a connected component of $\modk$. 
Any connected component of any $\modk$ has the form $\pi_{\R}(\Omega_0(\R))$ 
for some $\R$.
\item The map $\pi_{\R}$ is finite-to-one and almost everywhere locally 
bijective.
\item $\pi_{\R}(\U x)=\pi_{\R}(x)$.
\item The flow $P^t$ on $\Omega_0(\R)$ projects under $\pi_{\R}$ 
to the Teichm{\"u}ller flow $g_t$ on the corresponding connected 
component of $\modk$.
\item $(\pi_{\R})_*\mu_{\kappa}=\mu_{\R}$.
\item $m=2g-1+\sigma$.
\item  For any periodic orbit $\gamma$ of the flow $g_t$ on $\modk$, there exists a periodic orbit
${\tilde \gamma}$ of the flow $P^t$ on $\Omega_0(\R)$
such that $\pi_{\R}({\tilde \gamma})=\gamma.$
\end{enumerate}
\end{theorem}

\section{Appendix B: The uniform expansion property.}

This section is devoted to the uniform expansion property for the smooth 
invariant measure of the Teichm{\"u}ller flow on Veech's space of zippered rectangles.
I am deeply grateful to  Maryam Mirzakhani for explaining this property of the Teichm{\"u}ller flow to me.

First, Veech's coordinates on the space of zippered rectangles are modified, following \cite{bufetov}.  
Take a zippered rectangle $(\la, h, a, \pi)\in\Omega(\R)$, and introduce the vector
$\delta=(\delta_1, \dots, \delta_m)\in {\mathbb R}^m$ by the formula
$$
\delta_i=a_{i-1}-a_i, \ i=1, \dots, m
$$
(here we assume, as always, $a_0=a_{m+1}=0$).

We have then the following lemma from \cite{bufetov}.
\begin{proposition}
\label{deltazip}
 The data $(\la,\pi,\delta)$ determine the zippered rectangle $(\la,h, a,\pi)$ uniquely.
\end{proposition}

The parameters $h$ and $a$ of the zippered rectangle are expressed through the $\delta$
by the following formulas from \cite{bufetov}:
\begin{equation}
\label{hh}
h_{r}=-\sum_{i=1}^{r-1} \delta_i+\sum_{l=1}^{\pi(r)-1}\delta_{\pi^{-1}l}.
\end{equation}
\begin{equation}
\label{aa}
a_i=-\delta_1-\dots -\delta_{i-1}.
\end{equation}

In \cite{bufetov}, it is also proven that the relations (\ref{zipone})-(\ref{zipfour})
defining the zippered rectangle take the following equivalent form in the new coordinates:
\begin{equation}
\label{deltaone}
\delta_1+\dots +\delta_i\leq 0,\ \  i=1, \dots, m-1.
\end{equation}
\begin{equation}
\label{deltatwo}
\delta_{\pi^{-1}1}+\dots+\delta_{\pi^{-1}i}\geq 0, \ \  i=1, \dots, m-1.
\end{equation}
The parameter $a_m=-(\delta_1+\dots+\delta_m)$ can be both 
positive and negative. 

Introduce the following cone in ${\mathbb R}^m$:
$$
K_{\pi}=\{\delta=(\delta_1,\dots,  \delta_m): \delta_1+\dots +\delta_i\leq 0,
\delta_{\pi^{-1}1}+\dots+\delta_{\pi^{-1}i}\geq 0, i=1, \dots, m-1\},
$$
\begin{proposition} [\cite{bufetov}]
For $(\la,\pi)\in \Delta(\R)$ and an arbitrary $\delta\in K_{\pi}$
there exists a unique zippered rectangle $(\la,h,a,\pi)$ 
corresponding to the parameters $(\la,\pi,\delta)$.
\end{proposition}

In what follows, we shall simply refer to the zippered rectangle $(\la,\pi,\delta)$.

Denote by $Area(\la,\pi,\delta)$ the area of the zippered 
rectangle $(\la,\pi, \delta)$. We have from \cite{bufetov}:
\begin{equation}
\label{ziparea}Area(\la,\pi,\delta)=
\sum_{i=1}^m \delta_i (-\sum_{r=i+1}^m \la_r+\sum_{r=\pi(i)+1}^{m} 
\la_{\pi^{-1}r}).  
\end{equation}

Consider the set 
$$
{\cal V}(\R)=\{(\la,\pi,\delta):\pi\in\R, \la\in{\mathbb R}^m_+, \delta\in K(\pi)\}.
$$

In other words, ${\cal V}(\R)$ 
is the space of {\it all} possible zippered rectangles (i.e., not necessarily of 
those of area $1$).

The Teichm{\"u}ller flow  $P^t$ acts on ${\cal V}$ by the 
formula 
$$
P^t(\la,\pi,\delta)=(e^{t}\la, \pi, e^{-t}\delta).
$$

The map $\U$ acts on ${\cal V}(\R)$ by the formula

$$
{\U}(\la,\pi, \delta)=\begin{cases}
(A(\pi, b)^{-1}\la), b\pi, 
A(\pi, b)^{-1}\delta), 
&\text{if $\la\in\Delta_{\pi}^+$;}\\
(A(\pi, a)^{-1}\la), a\pi, 
A(\pi, a)^{-1}\delta), 
&\text{if $\la\in\Delta_{\pi}^-$.}
\end{cases}
$$

The volume form $Vol=d\la_1\dots d\la_md\delta_1\dots d\delta_m$ on ${\cal V}(\R)$
is preserved under the action of the flow $P^t$ and of the map $\U$.
Now consider the subset 
$$
{\cal V}^{(1)}(\R)=\{(\la,\pi,\delta): Area(\la,\pi,\delta)=1\},
$$ 
i.e., the subset of zippered rectangles of area $1$; naturally, this set is invariant under the 
Teichm{\"u}ller flow. 
Finally, to represent the moduli space of zippered rectangles in these new coordinates, 
take the quotient of ${\cal V}^{(1)}(\R)$ by the action of $\U$ and obtain the space 
 ${\cal V}^{(1)}_0(\R)$, which is isomorphic to $\Omega_0(\R)$, the linear isomorphism being given by the formulas 
(\ref{hh}), (\ref{aa}).

Now, given a Borel set $X\subset {\cal V}^{(1)}(\R)$, 
consider the set $K(X)\subset {\cal V}(\R)$, defined by the formula
\begin{equation}
\label{conus}
K(X)=\{t(\la,\delta, \pi), t\in [0,1], (\la,\delta, \pi)\in X\}.
\end{equation}
and introduce a measure ${\tilde \mu}_{\R}$ on ${\cal V}^{(1)}(\R)$ by the formula
$$
{\tilde \mu}_{\R}(X)=Vol(K(X)).
$$

Clearly, the measure ${\tilde \mu}_{\R}$ belongs to the Lebesgue measure class on $\vone$ 
and is invariant both under the action of flow $P^t$ and of the map $\U$. 
The corresponding quotient measure  ${\cal V}^{(1)}_0(\R)$ will be still denoted by  ${\tilde \mu}_{\R}$.
Under the isomorphism
of ${\cal V}^{(1)}_0(\R)$ and $\Omega_0(\R)$,
the measure ${\tilde \mu}_{\R}$ on  ${\cal V}^{(1)}_0(\R)$ is taken to 
(a constant multiple of) the measure $\mu_{\R}$.

Consider now two foliations on $\vone(\R)$:
$$
\F^-=\{x\in\vone(\R)| \ x=(\la,\delta,\pi): \la=const, \pi=const\},
$$
$$
\F^+=\{x\in\vone(\R)| \ x=(\la,\delta,\pi): \delta=const, \pi=const\}.
$$

Note that these foliations are invariant under the action of the map $\U$.

For a point $x\in\vone(\R)$, we denote by
$\F^+(x)$ the leaf of the foliation $\F^+$, passing
through $x$, by $\F^-(x)$ the leaf of the foliation $\F^-$, passing
through $x$.

Then, by  definition (\ref{conus})  of the measure
${\tilde \mu}_{\R}$, we have the following properties:
\begin{enumerate}
\item each leaf $\eta$ of the foliation $\F^+$
carries a globally defined conditional measure  ${\tilde \mu}^+_{\eta}$;
\item  each leaf $\xi$ of the foliation $\F^-$ carries
a globally defined conditional measure ${\tilde \mu}^-_{\xi}$;
\item $(P^t)_*({\tilde \mu}^+_{\eta})=exp(+mt){\tilde \mu}^+_{P^t\eta}$ for any $t>0$;
\item $(P^t)_*({\tilde \mu}^-_{\xi})=exp(-mt){\tilde \mu}^-_{P^t\xi}$ for any $t>0$.
\end{enumerate}

In other words, the measure ${\tilde \mu}_{\R}$ has the Margulis property 
\cite{margulis} of uniform expansion/contraction with respect to the pair 
of foliations $\F^+$, $\F^-$. In projection to the moduli space of abelian differentials, the 
foliations $\F^+$, $\F^-$ are taken to the stable and the unstable foliation of the 
Teichm{\"u}ller flow. The smooth measure on the moduli space thus also satisfies the Margulis property.

{\bf Acknowledgements.} 
I am deeply grateful to Maryam Mirzakhani for explaining to me 
the uniform expansion property for the Teichm{\"u}ller flow. 
I am deeply grateful to Artur Avila and Erwan Lanneau for many helpful suggestions. 
I am deeply grateful to Giovanni Forni, Vadim Kaimanovich, 
Grigorii Margulis and William Veech for useful discussions. 
While working on this project, I visited the Institute of Mathematics ``Guido Castelnuovo'' of 
the University of Rome ``La Sapienza''
at the invitation of Tullio Ceccherini-Silberstein, the Erwin
Schr{\"o}dinger Institute in Vienna at the invitation of Paul M{\"u}ller, 
and the Institut Henri Poincar{\'e} in the framework of the programme ``Time at work''. 
I am deeply grateful to these institutions for their hospitality.
In summer 2005 I was supported by the Clay Mathematics Institute's Liftoff Programme.


\begin{thebibliography}{9}

\bibitem{veech}
William Veech, Gauss measures for transformations on the space of interval exchange maps, Annals of Mathematics, 15(1982),
201--242.


\bibitem{veech2}
W.Veech, Interval exchange transformations.  J. Analyse Math.  33  (1978), 222--272.

\bibitem{veech3}
Veech, William A. Projective Swiss cheeses and uniquely ergodic interval exchange transformations.  Ergodic theory and
dynamical systems, I (College Park, Md., 1979--80),  pp. 113--193, 
Progr. Math., 10, Birkh{\"a}user, Boston, Mass., 1981.



\bibitem{zorich}

Anton Zorich, Finite Gauss measure on the space of interval exchange transformations. Lyapunov exponents. Ann. Inst. Fourier
(Grenoble) 46 (1996), no. 2, 325--370.


\bibitem{rauzy} G.Rauzy, {\'E}changes d'intervalles et transformations induites.
Acta Arith. 34,  (1979), no. 4, 315--328.

\bibitem{forni}
Giovanni Forni, Deviation of ergodic averages
for area-preserving flows on surfaces of higher genus.
Ann. of Math. (2),  155
(2002),  no. 1, 1--103.

\bibitem{masur} H. Masur, Interval exchange transformations and measured foliations. Ann. of Math. (2)  115  (1982), no. 1,
169--200.


\bibitem{sinailectures}
Ya.G. Sinai, Topics in ergodic theory, Princeton University Press, 1994.




\bibitem{Masur2} H. Masur, Interval exchange transformations and measured foliations, Ann. of Math., 115 (1982), 169--200.



\bibitem{bufetov} A. Bufetov, 
Decay of correlations for the Rauzy-Veech-Zorich induction map on 
the space of interval exchange transformations and the 
Central Limit Theorem for the Teichm{\"u}ller Geodesic Flow, 
preprint no. 1593 of the Erwin Schroedinger
Institute, Vienna, Jan. 2005; Ph.D. Thesis, 
Princeton, May 2005; arxiv-preprint, June 2005.

 \bibitem{veechteich} William Veech, The Teichm{\"u}ller geodesic flow, Annals of Mathematics, 1986.

\bibitem{margulis}
 G.A. Margulis, On some aspects of
the theory of Anosov flows, Ph.D. Thesis, 1970, Springer, 2003.

\bibitem{csf}

I.P. Cornfeld, Ya. G. Sinai, S.V. Fomin, Ergodic theory, 
``Nauka'', Moscow, 1980.

\bibitem{hubmas}
Hubbard, John; Masur, Howard, Quadratic differentials and foliations. Acta Math.  142  (1979), no. 3-4, 221--274.

\bibitem{kz1} M.Kontsevich, Lyapunov exponents and Hodge theory, ``Mathematical Beauty of Physics'', Saclay, 1996.

\bibitem{kontsevichzorich}
M.Kontsevich, A.Zorich, Connected components of the moduli spaces of Abelian differentials with prescribed singularities,
Inventiones mathematic\ae, 153(2003), no.3, 631-678.


\end{thebibliography}
\end{document}